\documentclass[12pt]{article} 
\usepackage{amsmath,graphicx,psfrag}
\def\R{I\!\!R}
\newcommand{\qed}{\nobreak \ifvmode \relax \else
\ifdim\lastskip<1.5em \hskip-\lastskip
\hskip1.5em plus0em minus0.5em \fi \nobreak
\vrule height0.75em width0.5em depth0.25em\fi}
\newcommand{\prll}{\, || \,}

\begin{document} 

\vspace*{-0.7in}

\begin{center}

{\Large{Eleven Euclidean Distances are Enough}} 

\vspace*{0.7cm}

{\large{Sujith Vijay}} 

\end{center}

\section*{Introduction}

The three distance theorem is a classic result in the study of
distributions modulo $1$, proved independently by several authors (see
\cite{Sos58} and \cite{Swi58}) in the 1950s in response to a conjecture of
Steinhaus. The theorem states that there are at most three distinct gaps
between consecutive elements in the set of fractional parts of the first
$n$ multiples of any real number $\alpha$. Formally, we have the
following: \\

{\bf {Theorem 1A}} Let $\alpha$ be any real number, and $n$ a positive
integer. Let $(a_1,a_2, \ldots, a_n)$ be the unique permutation of
$\{1,2,\ldots,n\}$ such that $$0 < \{a_1 \alpha \} < \{a_2 \alpha \} <
\cdots \{a_n \alpha \} < 1$$ Define $g_{\alpha}(0)= \{a_1 \alpha \}$ and
$g_{\alpha}(n)= 1 - \{a_n \alpha\}$. For $1 \le k \le n-1$, let
$g_{\alpha}(k) = \{a_{k+1} \alpha \} - \{a_k \alpha \}$. Define
$$S_{\alpha}(n)  = \{ g_{\alpha}(k): 0 \le k \le n \}$$ Then
$|S_{\alpha}(n)| \le 3$. \\

%\vspace*{-1.2in}

\begin{figure}[!h]
\begin{center}
\resizebox{6cm}{6cm}{\includegraphics{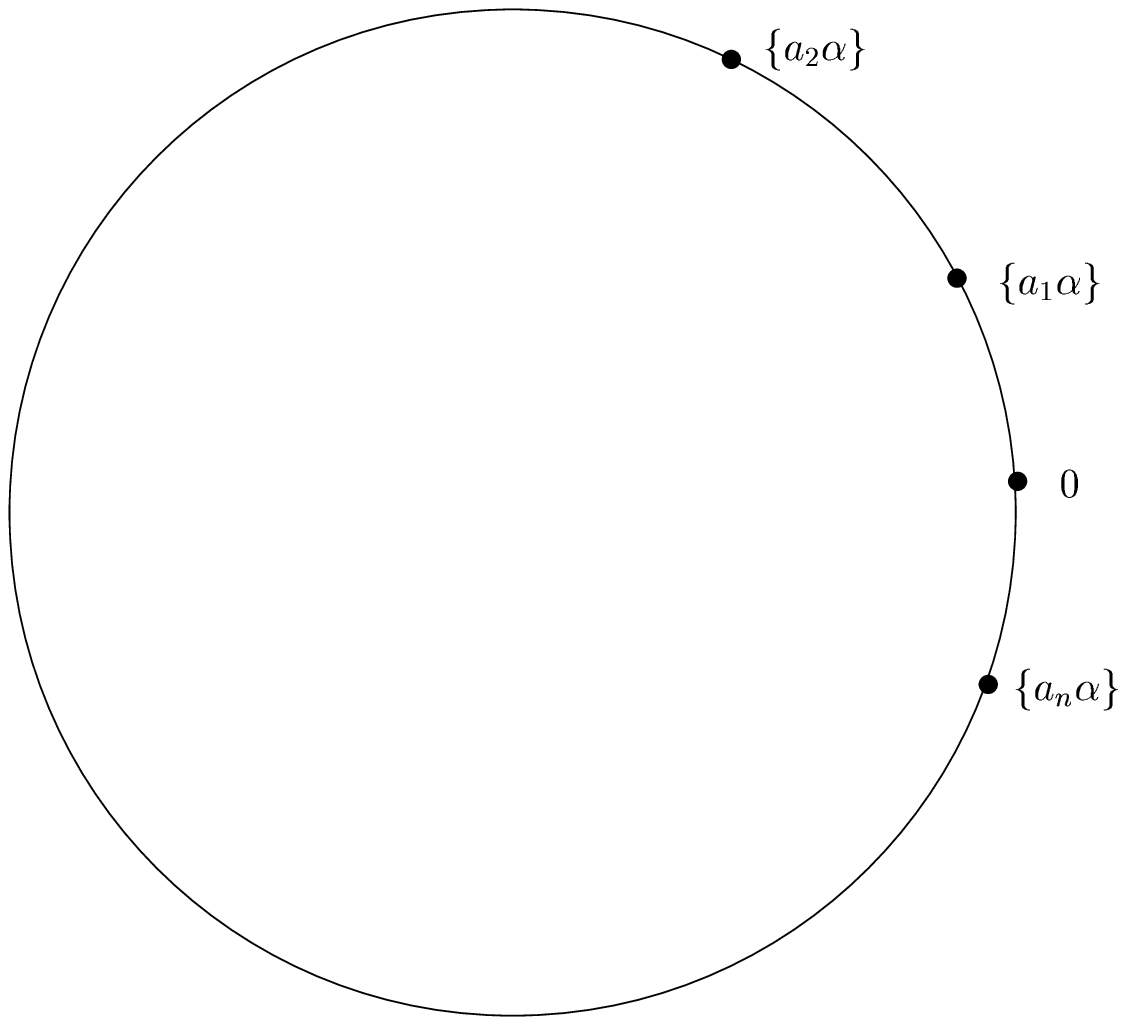}}
\end{center}
\end{figure}

\newpage

\section*{Generalisations}

Chung and Graham \cite{Chu76} generalised the three distance theorem as
follows: \\

{\bf {Theorem 2}} Let $\alpha, \lambda_1, \lambda_2, \ldots, \lambda_d$ be
real numbers, and let $n_1, n_2, \ldots, n_d$ be positive integers. For $1
\le i \le d, 1 \le k \le n$, let $a_{i,k} = \{ k \alpha + \lambda_i \}$,
where $\{x \}$ denotes the fractional part of $x$. Then there are at most
$3d$ distinct gaps between consecutive $a_{i,k}$. \\

Geelen and Simpson \cite{Gee93} established the following result, which 
was generalised by Chevallier \cite{Che00} to higher dimensions: \\

{\bf {Theorem 3}} Let $\alpha$ and $\beta$ be real numbers, and let $n_1$
and $n_2$ be positive integers. For $0 \le k_1 < n_1, 0 \le k_2 < n_2$,
let $a_{k_1,k_2} = \{ k_1 \alpha + k_2 \beta \}$. Then there are at most
$n_1 + 3$ distinct gaps between consecutive $a_{k_1,k_2}$. \\

Chevallier \cite{Che97} also obtained the following higher-dimensional
analogue of the three-distance theorem for a certain subsequence of
natural numbers. \\

{\bf {Theorem 4}} Let $N$ be a best simultaneous approximation denominator 
with respect to the Euclidean norm of the $d$-tuple $(\alpha_1, \alpha_2, 
\ldots, \alpha_d)$. Then there is a norm on $R^{d}$ such that the Voronoi 
diagrams of the first $N$ points of the sequence $(\{k \alpha_1 \}, \{k 
\alpha_2 \}, \ldots, \{k \alpha_d \})$ with respect to this norm are of at 
most $C_d$ different forms, where $C_d$ is a constant that depends only on 
the dimension $d$. 

\section*{A New Formulation}

The purpose of this article is to show that the central tenet of the
three-distance theorem, namely the finiteness of the set of minimal
distances, can be generalised to higher dimensions under a suitable
interpretation. We begin by rephrasing the theorem in a form that lends
itself to the generalisation we seek. \\

We think of the three distance theorem as a statement about champions in a
tournament. The players in the tournament are edges connecting $\{j
\alpha\}$ and $\{k \alpha\}, 1 \le j < k \le n$, two edges play each other
if and only if they overlap, and an edge loses only against edges of
shorter length that it plays against. Defeated edges are allowed to play
(and defeat) other overlapping edges. According to the three distance
theorem, there are at most three distinct values for the lengths of
undefeated edges. Thus the theorem can be restated as follows: \\

{\bf {Theorem 1B}} Let $\alpha$ be any real number, and $n$ a positive
integer. Define $d_{\alpha}(j,k)=||\{k \alpha \} - \{j \alpha \}||$. Let
$I_{j,k}$ be the ``geodesic" joining $\{j \alpha \}$ with $\{k \alpha \}$,
i.e., if $m_{j,k}= \min(\{j \alpha \}, \{k \alpha \}) \mbox{ and }
M_{j,k}= \max(\{j \alpha \}, \{k \alpha \})$, we define $$I_{j,k} =
\left\{ \begin{array}{ll} [m_{j,k},M_{j,k}) & \quad \mbox{if $M_{j,k} -
m_{j,k} \le 1/2$} \\ {[0,m_{j,k}) \cup [M_{j,k},1)} & \quad
\mbox{otherwise} \\ \end{array} \right.$$ Let $S_{\alpha}(n)= \{
d_{\alpha}(j,k):  d_{\alpha,\beta}(p,q) < d_{\alpha,\beta}(j,k)  
\Rightarrow I_{p,q} \cap I_{j,k} = \emptyset \}$. Then $|S_{\alpha}(n)|
\le 3$. \\

We first prove a two-dimensional version of this theorem. We show that if
the players are edges connecting $(\{j \alpha \}, \{j \beta\})$ and $(\{k
\alpha \}, \{ k \beta \})$ and two edges play each other if and only if
their {\em {projections along either co-ordinate axis}} overlap, there are
at most $11$ distinct values for the lengths of undefeated edges.
Numerical evidence suggests that the true value could be as small as $3$. \\

{\bf {Theorem 1}} Let $\alpha$ and $\beta$ be real numbers, and let $n$ be
a positive integer.  Define $d_{\alpha,\beta}(j,k)=\sqrt{||(k-j)
\alpha||^2 + ||(k-j) \beta ||^2}$. Let $I^1_{j,k}$ and $I^2_{j,k}$ be the
geodesics joining $\{j \alpha \}$ with $\{ k \alpha \}$ and $\{j \beta \}$
with $\{k \beta \}$ respectively. Define $$S_{\alpha,\beta}(n)=
\{d_{\alpha,\beta}(j,k):  d_{\alpha,\beta}(p,q) < d_{\alpha,\beta}(j,k)
\Rightarrow I^1_{p,q} \cap I^1_{j,k} = I^2_{p,q} \cap I^2_{j,k} =
\emptyset \}$$ Then $|S_{\alpha,\beta}(n)| \le 11$. \\

{\bf {Proof}} We begin by classifying the denominators of simultaneous
rational approximations to $(\alpha, \beta)$. Let $[[x]] = \{x\} - 1/2$.  
We say that $q$ is a denominator of type $(+,-)$ if $[[q \alpha]] \ge 0$
and $[[q \beta]] < 0$. Denominators of type $(-,+), (+,+)$ and $(-,-)$ are
defined analogously. Types $(+,+)$ and $(-,-)$ are said to be opposites to
each other, as are types $(+,-)$ and $(-,+)$. We write $q_1 \prll q_2$ if
$q_1$ and $q_2$ are of the same type, $q_1 \perp q_2$ if they are of
opposite type, and $q_1 \sim q_2$ if they are {\em {not}} of opposite
type. \\

We define the {\em {length}} and the {\em {angle}} of an integer $q$ with 
respect to $\alpha$ and $\beta$ as $\ell(q)=d_{\alpha,\beta}(0,q)$, and
$\theta(q)= \tan^{-1}([[q \beta]] / [[q \alpha]])$ respectively. Let $Q_1$
denote the least integer with the property that $\ell(Q_1) \le \ell(q)$
for all $q, \, 1 \le q \le n/2$.  For $n/2 < q \le n$, we say that $q$ is
{\em {primary}} if $\ell(q) < \ell(Q_1)$. \\

{\bf {Lemma 1}} If $q$ is primary, $\ell(q)$ can take at most five
distinct values. \\

{\bf {Proof}} Consider four quarter-circles of radius $R=\ell(Q_1)$ centred 
at the four corners of the unit square. Suppose there exist $q_i, 1 \le i \le 7$ with $n/2 < q_1 <
q_2 < \cdots < q_7 \le n$ and $\ell(q_i) < R$. Then there must be a pair
$(i,j), 1 \le i < j \le 7$ such that $\theta \doteq |\theta(q_j) - \theta(q_i)| < \pi/3$. But then we have $\ell(q_j - q_i) < R$, yielding a contradiction,
since $1 \le q_j - q_i < n/2$. \\

\begin{figure}[!h]
\begin{center}
\resizebox{6cm}{6cm}{\includegraphics{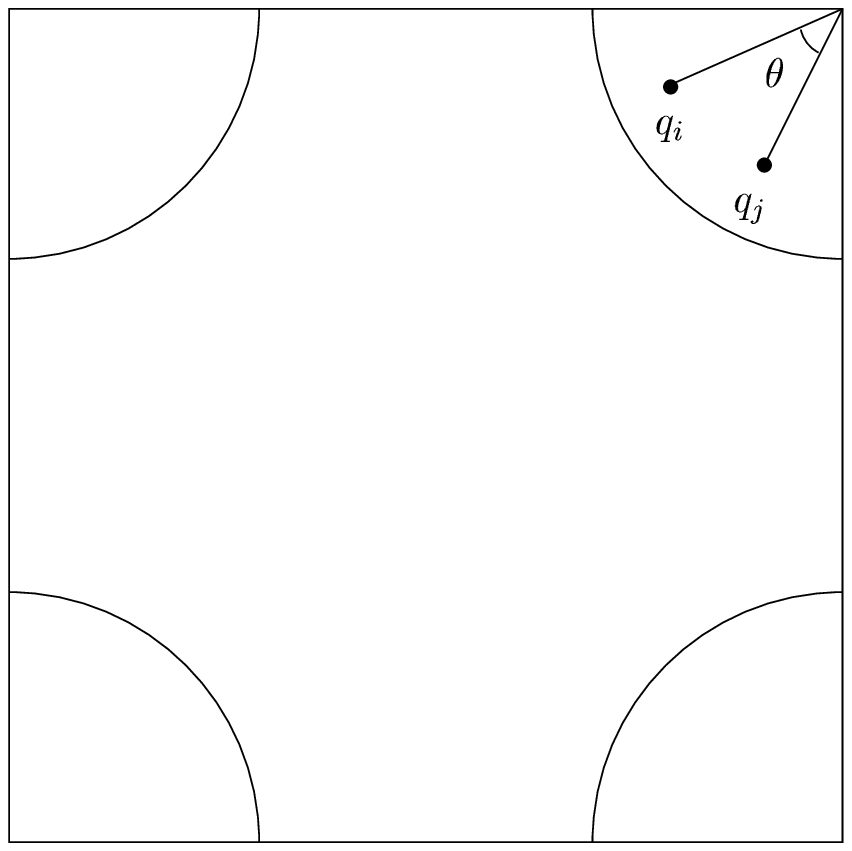}}
\end{center}
\end{figure}

Furthermore, the only way to have six primary $q_i$ avoiding $\ell(q_j-q_i)
 < R$ is to arrange them along the vertices of a regular hexagon, leading to 
identical values of $\ell(q_i)$. It follows that $\ell(q)$ can take at 
most five distinct values if $q$ is primary. \qed \\

Consider the line $L_{a,b}$ joining $(\{a \alpha \}, \{a \beta \})$ and
$(\{b \alpha \}, \{b \beta \})$, with $1 \le a < b \le n$. Let
$q^{*}=b-a$. As we have seen, if $q^{*}$ is primary, there are
only five possible values for $\ell(q^{*})$. Suppose $q^{*}$ is not
primary. We consider two cases. \\

{\small{\bf {CASE 1: $\mathbf{q^{*} \sim Q_1}$}}} Note that one of
$L_{a,a+Q_1}$ or $L_{b-Q_1,b}$ will be admissible, and will defeat
$L_{a,b}$. \\

{\small{\bf {CASE 2: $\mathbf{q^{*} \perp Q_1}$}}} Define $Q_1^{\perp} =
\{q: 1 \le q \le n-Q_1, q \not \! \! \! \prll Q_1 \}$. Note that
$Q_1^{\perp}$ is non-empty, since it contains $q^{*}$. Let $Q_2$
be the least integer in $Q_1^{\perp}$ such that $\ell(q)  \le \ell(Q_2)$
for all $q \in Q_1^{\perp}$. We first prove the following lemma. \\

{\bf {Lemma 2}} There is at most one $q < Q_1$ satisfying $\ell(q) <
\ell(Q_2)$. \\

{\bf {Proof}} Suppose there are at least two such $q$. Let $q_0$ and
$q'_0$ be the two smallest, with $q_0 < q'_0$. By the definition of $Q_2$,
$q_0 \prll q'_0 \prll Q_1$. 

Let $x_1 = ||Q_1 \alpha||, x_2 = ||Q_1 \beta||, y_1 = ||q_0 \alpha$ and
$y_2 = ||q_0 \beta||$. Since $\ell(Q_1)  \le \ell(q_0)$, we have
$$\ell(Q_1-q_0) = \sqrt{(x_1-y_1)^2 + (x_2 - y_2)^2} \le \sqrt{2(y_1^2 +
y_2^2)} = \sqrt{2} \ell(q_0) < \sqrt{2} \ell(Q_2)$$

Since $Q_1 - q_0 \not \! \! \! \prll Q_1$, it follows from the definition 
of $Q_2$ that $\ell(Q_1 - q_0) \ge \ell(Q_2) > \ell(q_0)$.  Similarly, 
$\ell(Q_1 - q'_0) \ge \ell(Q_2) > \ell(q'_0)$. This requires 
$\theta_0 \doteq |\theta(Q_1) - \theta(q_0) | \ge \pi/3$ and $\theta'_0 
\doteq |\theta(Q_1) - \theta(q'_0) | \ge \pi/3$. \\

\begin{figure}[!h]
\begin{center}
\resizebox{14cm}{2cm}{\includegraphics{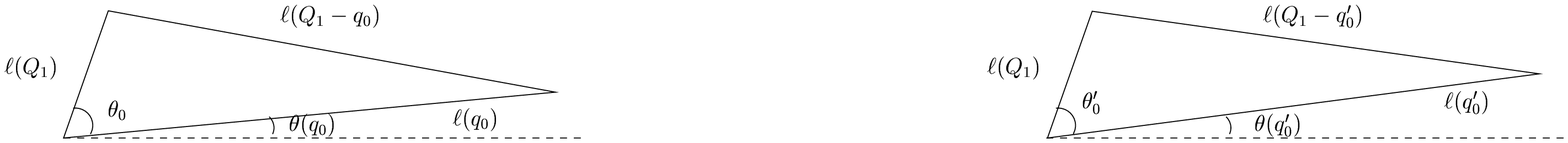}}
\end{center}
\end{figure}

It follows that $|\theta(q'_0) - \theta(q_0)| \le \pi/6$, so
$\ell(q'_0-q_0) < \ell(Q_2)$. Since $\ell(q_0) < \ell(Q_2)$, we must have
$q'_0 = 2 q_0$. Since $q_0 \prll q'_0$, we have, $$\ell(q'_0)=2 \ell(q_0)
\ge \sqrt{2} \ell(Q_1-q_0) \ge \sqrt{2} \ell(Q_2)$$ contradicting our
assumption about $q'_0$. It follows that there is at most one $q < Q_1$
satisfying $\ell(q) < \ell(Q_2)$. \qed \\

For $n-Q_1 < q \le n$ and $q \perp Q_1$, we say that $q$ is {\em
{secondary}} if $\ell(q) < \ell(Q_2)$. If $q^{*}$ is not secondary, one of
$L_{a,a+Q_2}$ or $L_{a-Q_1,a}$ will be admissible, and will defeat
$L_{a,b}$. We claim that if $q^{*}$ is secondary, there are at most four
distinct values that $q^{*}$ can take. \\

Suppose not. Let $q_1 > q_2 > \cdots > q_5 > n-Q_1$, with $\ell(q_i) <
\ell(Q_2)$ and $q_i \perp Q_1$. Note that if $|\theta(q_i) - \theta(q_j)|
\le \pi/3$, we have $\ell(q_i - q_j) \le \max(\ell(q_i), \ell(q_j)) <
\ell(Q_2)$. \\

But given $q_1 \prll q_2 \prll \ldots \prll q_5$, it is easy to see that
there exist $(i_1,j_1)$ and $(i_2,j_2)$ satisfying $q_{i_1} - q_{j_1} \neq
q_{i_2} - q_{j_2}$ and $$\max(|\theta(q_{i_1}) - \theta(q_{j_1})|,
|\theta(q_{i_2})  - \theta(q_{j_2})|) \le \pi/4 < \pi/3 $$

Thus we have $\ell(q_{i_1} - q_{j_1}) < \ell(Q_2)$ and $\ell(q_{i_2} -
q_{j_2})  < \ell(Q_2)$, contradicting Lemma 2. So there are at most four
distinct values that $\ell(q^{*})$ can take if $q^{*} \perp Q_1$ and
$q^{*}$ is secondary. It follows that at most eleven distinct gaps
survive, completing the proof of the theorem. \qed

\section*{Higher Dimensions}

For higher dimensions, the above argument can be adapted to obtain similar
results. We prove the following theorem which implies, in particular, that
there are at most 290 distances in three dimensions. \\

{\bf {Theorem 2}} Let $\mathbf{\alpha} \doteq (\alpha_1, \alpha_2, \ldots,
\alpha_m) \in {\R}^m$, and let $n$ be a positive integer.  Define
$$d_{\mathbf{\alpha}}(j,k)=\sqrt{\sum_{i=1}^n ||(k-j) \alpha_i ||^2}$$ For
$1 \le r \le m$, let $I^r_{j,k}$ denote the geodesic joining $\{j \alpha_r
\}$ with $\{ k \alpha_r \}$. Define $$S_{\mathbf{\alpha}}(n)=
\{d_{\mathbf{\alpha}}(j,k):  d_{\mathbf{\alpha}}(p,q) <
d_{\mathbf{\alpha}}(j,k) \Rightarrow I^r_{p,q} \cap I^r_{j,k} = \emptyset
\mbox{ for all } r. \} $$ Then $$|S_{\mathbf{\alpha}}(n)| \le ({\lceil
\sqrt{m} \, \rceil}^m)({\lceil \sqrt{m} \, \rceil}^m + 2^m + 1) + 2$$

{\bf {Proof}} Let $[[x]] = \{x\} - 1/2$. As in the proof of Theorem 1, we
assign, to each denominator $q$ an $m$-tuple of signs. The $i^{th}$ sign
is positive if $[[q \alpha_i]] \ge 0$ and negative otherwise. \\

The {\em {length}} of an integer $q$ with respect to $\mathbf{\alpha}$ is defined as $\ell(q)=d_{\mathbf{\alpha}}(0,q)$. Let $Q_1$ denote the least integer with the 
property that $\ell(Q_1) \le \ell(q)$ for all $q, \, 1 \le q \le n/2$.  For 
$n/2 < q \le n$, we say that $q$ is {\em {primary}} if $\ell(q) < \ell(Q_1)$. \\

{\bf {Lemma 3}} There are at most ${(2 \lceil \sqrt{m} \, \rceil)}^m$
primary denominators in $\R^m$ for any given $\mathbf{\alpha}$. \\

{\bf {Proof}} If the number of distinct values $q$ satisfying $\ell(q) <
\ell(Q_1)$ exceeds ${(2 \lceil \sqrt{m} \, \rceil)}^m$, at least $1 +
{\lceil \sqrt{m} \, \rceil}^m$ of these values must be of the same type.
By pigeonhole principle, there exists $q_1$ and $q_2$ with $||(q_2 - q_1)
\alpha_i|| < \ell(Q_1) /\sqrt{m}$ for all $i, 1 \le i \le m$. It follows
that $\ell(q_2-q_1) < \ell(Q_1)$. But $q_2 - q_1 < n/2$, contradicting the
definition of $Q_1$. Thus at most ${(2 \lceil \sqrt{m} \, \rceil)}^m$
denominators can be primary. \qed \\

Consider the line $L_{a,b}$ joining $(\{a \alpha \}, \{a \beta \})$ and
$(\{b \alpha \}, \{b \beta \})$, with $1 \le a < b \le n$. Let
$q^{*}=b-a$. As we have seen, if $q^{*}$ is primary, there are
only five possible values for $\ell(q^{*})$. Suppose $q^{*}$ is not
primary. We consider two cases. \\

{\small{\bf {CASE 1: $\mathbf{q^{*} \sim Q_1}$}}} Note that one of
$L_{a,a+Q_1}$ or $L_{b-Q_1,b}$ will be admissible, and will defeat
$L_{a,b}$. \\

{\small{\bf {CASE 2: $\mathbf{q^{*} \perp Q_1}$}}} As in the planar case,
define $Q_1^{\perp} = \{q: 1 \le q \le n-Q_1, q \not \! \! \! \prll Q_1
\}$, and let $Q_2$ be the least integer in $Q_1^{\perp}$ such
that $\ell(q)  \le \ell(Q_2)$ for all $q \in Q_1^{\perp}$. We prove
the following analogue of Lemma 2. \\

{\bf {Lemma 4}} There are at most ${\lceil \sqrt{2m} \, \rceil}^m$ values
of $q < Q_1$ satisfying $\ell(q) < \ell(Q_2)$. \\

{\bf {Proof}} Suppose there are more. Then there exist $q_0$ and $q'_0$
with $||(q'_0 - q_0) \alpha_i|| < \ell(Q_2) / \sqrt{2d}$ for all $i, 1 \le
i \le m$. Thus $\ell(Q_2) > \sqrt{2} \ell(q''_0)$ where $q''_0 = q'_0 -
q_0$. Note that $Q_1 - q''_0 \not \! \! \! \prll Q_1$.  Therefore
$\ell(Q_2) < \ell(Q_1 - q''_0$. We will now deduce a contradiction by
showing that $\ell(Q_1 - q''_0) \le \sqrt{2} \ell(q''_0)$. \\

Let $x_r= ||Q_1 \alpha_r|| and y_r = ||q''_0 \alpha_r||$. Since $\ell(Q_1)  
\le \ell(q''_0)$, we have $$\ell(Q_1-q_0) = \sqrt{\sum_{r=1}^m
(x_r-y_r)^2} \le \sqrt{2 \sum_{r=1}^m y_r^2} = \sqrt{2}
\ell(q''_0)$$ yielding the desired contradiction. It follows that there
are at most ${\lceil \sqrt{2m} \, \rceil}^m$ values of $q < Q_1$
satisfying $\ell(q) < \ell(Q_2)$ \qed \\

For $n-Q_1 < q \le n$ and $q \perp Q_1$, we say that $q$ is {\em
{secondary}} if $\ell(q) < \ell(Q_2)$. If $q^{*}$ is not secondary, one of
$L_{a,a+Q_2}$ or $L_{a-Q_1,a}$ will be admissible, and will defeat
$L_{a,b}$. We claim that if $q^{*}$ is secondary, there are at most
$({\lceil \sqrt{m} \, \rceil}^m) ({\lceil \sqrt{2m} \, \rceil}^m + 1)$
distinct values that $q^{*}$ can take.\\

Suppose not. Let $k=({\lceil \sqrt{m} \, \rceil}^m) ({\lceil \sqrt{2m} \,
\rceil}^m + 1) + 1$, and let $$q_1 > q_2 > \cdots > q_k > n-Q_1$$ with
$\ell(q_i) < \ell(Q_2)$. Then there exist ${\lceil \sqrt{2m} \, \rceil}^m
+ 1$ distinct denominators $q < Q_1$ with $||q \alpha_r|| < \ell(Q_2) /
\sqrt{m}$ for all $r, 1 \le r \le m$, thus satisfying $\ell(q)  <
\ell(Q_2)$ and contradicting Lemma 4.  So there are at most $({\lceil
\sqrt{m} \, \rceil}^m) ({\lceil \sqrt{2m} \, \rceil}^m + 1)$ distinct
values that $\ell(q^{*})$ can take if $q^{*} \perp Q_1$ and $q^{*}$ is
secondary. Accounting for $Q_1$, $Q_2$ and primary denominators, we obtain
the statement of the theorem. \qed

\end{document}